\documentclass[12pt]{article}

\newtheorem{theorem}{Theorem}

\newtheorem{example}[theorem]{Example}

\newtheorem{remark}[theorem]{Remark}

\newenvironment{proof}[1][Proof]{\textbf{#1.} }{\ \rule{0.5em}{0.5em}}
\newdimen\dummy
\dummy=\oddsidemargin
\date{}
\begin{document}

\title{A characterization of spaces with discrete topological fundamental group}
\author{Paul Fabel \\
Department of Mathematics \& Statistics\\
Mississippi State University}
\maketitle

\begin{abstract}
The fundamental group of a locally path connected metric space inherits the
discrete topology in a natural way if and only if $X$ is semilocally simply
connected. We also provide a counterexample to a similar theorem in the
literature.
\end{abstract}

\section{Introduction}

In general $\pi _{1}(X,p),$ the based fundamental group of a space $X,$
admits a canonical topology and becomes the \textit{topological fundamental
group.} To date significant progress has been made in the investigation of
the topological fundamental group of certain non semilocally simply
connected spaces such as the Hawaiian earring and the harmonic archipelago.
For example, the Hawaiian earring group is not a Baire space \cite{fab3} and
also fails to embed topologically in the inverse limit of free groups \cite
{fab25} despite injectivity of the canonical homomorphism \cite{mor}. The
harmonic archipelago is a compact metric space whose fundamental group is
uncountable but has only two open sets \cite{fab2}. In particular neither of
these groups has the discrete topology.

When does $\pi _{1}(X,p)$ have the discrete topology? The fundamental groups
of locally contractible spaces such as $n-$manifolds have the discrete
topology. Indeed, a published theorem of another author (Theorem 5.1 \cite
{Biss}) indicates that $\pi _{1}(X,p)$ has the discrete topology precisely
when $X$ is semilocally simply connected.

Unfortunately Theorem 5.1 is false. This paper contains a counterexample, an
alternate version of the theorem, and some remarks and examples that lend
further perspective on the results in question.

\section{Definitions}

Suppose $X$ is a metrizable space and $p\in X.$ Let $C_{p}(X)=\{f:[0,1]%
\rightarrow X$ such that $f$ is continuous and $f(0)=f(1)=p\}.$ Endow $%
C_{p}(X)$ with the topology of uniform convergence. The \textbf{topological
fundamental group} $\pi _{1}(X,p)$ is the set of path components of $%
C_{p}(X) $ topologized with the quotient topology under the canonical
surjection $q:C_{p}(X)\rightarrow \pi _{1}(X,p)$ satisfying $q(f)=q(g)$ if
and only of $f $ and $g$ belong to the same path component of $C_{p}(X).$

Thus a set $U\subset \pi _{1}(X,p)$ is open in $\pi _{1}(X,p)$ if and only
if $q^{-1}(U)$ is open in $\pi _{1}(X,p).$

A space $X$ is \textbf{semilocally simply connected} at $p$ if $p\in X$ and
there exists an open set $U$ such that $p\in U$ and $j_{U}:U\rightarrow X$
induces the trivial homomorphism $j_{U}^{\ast }:\pi _{1}(U,p)\rightarrow \pi
_{1}(X,p).$ Roughly speaking, this means small loops in $X$ are inessential,
but might bound large disks.

A topological space $X$ is \textbf{discrete} if every subset is both open
and closed. The above definitions are consistent with those found in Munkres 
\cite{Munkres}.

\section{A counterexample}

Looking beyond the notational ambiguity in case $X$ is not path connected,
it is claimed falsely in Theorem 5.1 \cite{Biss} that $\pi _{1}(X)$ is
discrete if and only if the topological space $X$ is semilocally simply
connected.

Here is a counterexample to Theorem 5.1. 

\begin{example}
\label{e1}Let $X$ denote the following subset of the plane: 
\[
X=(\{0,1,\frac{1}{2},\frac{1}{3},...\}\times \lbrack 0,1])\cup ([0,1]\times
\{0\})\cup ([0,1]\times \{1\}).
\]
Let $p=(0,0).$ Consider the inessential map $f$ which starts at $(0,0),$
goes straight up to $(0,1)$ and then comes back down again. Notice we may
construct maps $f_{n}\in C_{p}(X)$ such that $f_{n}$ is inessential but $%
f_{n}\rightarrow f$ uniformly. (Let $f_{n}$ trace a rectangle with sides $%
\{0\}\times \lbrack 0,1]$ and $\{\frac{1}{n}\}\times \lbrack 0,1]$). This
shows that the path component of the constant map is not open in $C_{p}(X).$
Consequently $\pi _{1}(X,p)$ does not have the discrete topology. On the
other hand $X$ is semilocally simply connected since small loops in $X$ are
null homotopic in $X.$
\end{example}

\section{A characterization of discrete topological fundamental groups}

Theorem 5.1 \cite{Biss} can be repaired by making the added assumption that $%
X$ is locally path connected. 

\begin{theorem}
\label{main}Suppose $X$ is a locally path connected metrizable space and $%
p\in X.$ Then $\pi _{1}(X,p)$ is discrete if and only if $X$ is semilocally
simply connected at $p.$
\end{theorem}

\begin{proof}
Let $F:C_{p}(X)\rightarrow \pi _{1}(X,p)$ denote the canonical quotient map.
Let $P\in C_{p}(X)$ denote the constant map. Let $[P]$ denote the path
component of $P$ in $C_{p}(X).$ Let $d$ be any metric on $X.$

Suppose $\pi _{1}(X,p)$ is discrete. Then the trivial element of $\pi
_{1}(X,p)$ is an open subset. Hence, since $F$ is continuous, $[P]$ is open
in $C_{p}(X).$ Suppose, in order to obtain a contradiction, that $\pi
_{1}(X,p)$ is not semilocally simply connected at $p.$ Then there exists a
sequence of maps $\alpha _{n}\rightarrow P$ such that $\alpha _{n}\in
C_{p}(X)$ and $\alpha _{n}$ is essential in $X$ for each $n.$ On the other
hand, since $[P]$ is open in $C_{p}(X)$ we can be sure that eventually $%
\alpha _{n}$ is inessential in $X.$ 

Conversely, suppose $X$ is semilocally simply connected at $p.$ To prove $%
\pi _{1}(X,p)$ is discrete we must show that each one point subset is open.
Since $\pi _{1}(X,p)$ is a topological group (Proposition 3.1 \cite{Biss}) $%
\pi _{1}(X,p)$ is homogeneous. Thus it suffices to prove that the trivial
element of $\pi _{1}(X,p)$ is an open subset. Since $F$ is a quotient map it
suffices to prove $[P]$ is open in $C_{p}(X).$ Suppose $f\in \lbrack p]$ and
suppose $f_{n}\in C_{p}(X)$ and suppose $f_{n}\rightarrow f$ uniformly. We
must prove that for large $n$ $f$ and $f_{n}$ are path homotopic.

Observation1. Since $X$ is locally path connected and since $im(f)$ is
compact, for each $\varepsilon >0$ there exists $\delta >0$ such that if $%
x\in im(f)$ and $d(x,y)<\delta $ then there exists a path $\alpha _{xy}$
from $x$ to $y$ such that diam$(\alpha _{xy})<\varepsilon .$ 

Observation 2. Since $X$ is semilocally simply connected and since $im(f)$
is compact, there exists $\varepsilon >0$ such that if $x\in im(f)$ and if $%
\alpha _{xx}$ is any path from $x$ to $x$ such that diam$(\alpha
_{xx})<\varepsilon $ and then $\alpha _{xx}$ is null homotopic. 

Observation 3. The set $f\cup \{f_{1},f_{2},..\}$ is an equicontinuous
collection of maps.

Combining observations 1,2, and 3, for sufficiently large $n$ we can
manufacture a homotopy from $f_{n}$ to $f$ described roughly as follows.
Observation 1 says we can connect $f(\frac{i}{2^{n}})$ and $f_{n}(\frac{i}{%
2^{n}})$ by a small path and we can connect $f(\frac{i+1}{2^{n}})$ and $%
f_{n}(\frac{i+1}{2^{n}})$ by a small path. Observation 3 says the
restriction of $f$ and $f_{n}$ to $[\frac{i}{2^{n}},\frac{i+1}{2^{n}}]$ is
small. Observation 2 says the boundary of the rectangle with corners  $f(%
\frac{i}{2^{n}})$, $f_{n}(\frac{i}{2^{n}}),f(\frac{i+1}{2^{n}})$ and $f_{n}(%
\frac{i+1}{2^{n}})$ determines an inessential loop. Thus we can fill in the
homotopy across $[\frac{i}{2^{n}},\frac{i+1}{2^{n}}]\times \lbrack 0,1].$
\end{proof}

\begin{remark}
Local path connectivity is only used in the 2nd part of the proof of Theorem 
\ref{main}. Thus all path connected spaces with discrete topological
fundamental group must be semilocally simply connected.
\end{remark}

\begin{remark}
There exist spaces which are not locally path connected but which have
discrete fundamental group. For example the cone over $X$ in example \ref{e1}%
. (i.e. let $Y$ be quotient of the space $X\times \lbrack 0,1]$ with $%
X\times \{1\}$ considered as a point. Note $Y$ is simply connected.)
\end{remark}

\begin{remark}
Theorem 5.1 of \cite{Biss} cannot be repaired by replacing the notion of
`discrete' with `totally disconnected'. The space $X=S^{1}\times S^{1}..$
has totally disconnected fundamental group but $X$ is not semilocally simply
connected.
\end{remark}

\end{document}